\documentclass[12pt]{article}

\usepackage{amsmath,amssymb}
\begin{document}

\pagestyle{myheadings} \markright{LEFSCHETZ FORMULAE FOR p-ADIC GROUPS}

\title{Lefschetz formulae for $p$-adic groups}
\author{Anton Deitmar}

\date{}
\maketitle

$$ $$

\tableofcontents

\def \1{{\bf 1}}
\def \a{{{\mathfrak a}}}
\def \ad{{\rm ad}}
\def \al{\alpha}
\def \ar{{\alpha_r}}
\def \A{{\mathbb A}}
\def \Ad{{\rm Ad}}
\def \Aut{{\rm Aut}}
\def \b{{{\mathfrak b}}}
\def \bs{\backslash}
\def \B{{\cal B}}
\def \c{{\mathfrak c}}
\def \cent{{\rm cent}}
\def \C{{\mathbb C}}
\def \CA{{\cal A}}
\def \CB{{\cal B}}
\def \CC{{\cal C}}
\def \CD{{\cal D}}
\def \CE{{\cal E}}
\def \CF{{\cal F}}
\def \CG{{\cal G}}
\def \CH{{\cal H}}
\def \CHC{{\cal HC}}
\def \CL{{\cal L}}
\def \CM{{\cal M}}
\def \CN{{\cal N}}
\def \CP{{\cal P}}
\def \CQ{{\cal Q}}
\def \CO{{\cal O}}
\def \CS{{\cal S}}
\def \CT{{\cal T}}
\def \CV{{\cal V}}
\def \CW{{\cal W}}
\def \det{{\rm det}}
\def \df{\ \begin{array}{c} _{\rm def}\\ ^{\displaystyle =}\end{array}\ }
\def \diag{{\rm diag}}
\def \dist{{\rm dist}}
\def \End{{\rm End}}
\def \eps{\varepsilon}
\def \eqn{\begin{eqnarray*}}
\def \endeqn{\end{eqnarray*}}
\def \F{{\mathbb F}}
\def \Fx{{\mathfrak x}}
\def \FX{{\mathfrak X}}
\def \g{{{\mathfrak g}}}
\def \ga{\gamma}
\def \Ga{\Gamma}
\def \GL{{\rm GL}}
\def \h{{{\mathfrak h}}}
\def \Hom{{\rm Hom}}
\def \im{{\rm im}}
\def \Im{{\rm Im}}
\def \Ind{{\rm Ind}}
\def \k{{{\mathfrak k}}}
\def \K{{\cal K}}
\def \l{{\mathfrak l}}
\def \la{\lambda}
\def \lap{\triangle}
\def \li{{\rm li}}
\def \Lie{{\rm Lie}}
\def \m{{{\mathfrak m}}}
\def \mod{{\rm mod}}
\def \n{{{\mathfrak n}}}
\def \name{\bf}
\def \Mat{{\rm Mat}}
\def \N{\mathbb N}
\def \o{{\mathfrak o}}
\def \ord{{\rm ord}}
\def \O{{\cal O}}
\def \p{{{\mathfrak p}}}
\def \ph{\varphi}
\def \prf{\noindent{\bf Proof: }}
\def \Per{{\rm Per}}
\def \PGL{{\rm PGL}}
\def \q{{\mathfrak q}}
\def \qed{\ifmmode\eqno $\square$\else\noproof\vskip 12pt plus 3pt minus 9pt \fi}
 \def\noproof{{\unskip\nobreak\hfill\penalty50\hskip2em\hbox{}%
     \nobreak\hfill $\square$\parfillskip=0pt%
     \finalhyphendemerits=0\par}}
\def \Q{\mathbb Q}
\def \res{{\rm res}}
\def \R{{\mathbb R}}
\def \Re{{\rm Re \hspace{1pt}}}
\def \r{{\mathfrak r}}
\def \ra{\rightarrow}
\def \rank{{\rm rank}}
\def \sign{{\rm sign}}
\def \supp{{\rm supp}}
\def \SL{{\rm SL}}
\def \Spin{{\rm Spin}}
\def \SU{{\rm SU}}
\def \t{{{\mathfrak t}}}
\def \T{{\mathbb T}}
\def \tr{{\hspace{1pt}\rm tr\hspace{1pt}}}
\def \vol{{\rm vol}}
\def \z{\zeta}
\def \Z{\mathbb Z}
\def \={\ =\ }

\newcommand{\frack}[2]{\genfrac{}{}{0pt}{}{#1}{#2}}
\newcommand{\rez}[1]{\frac{1}{#1}}
\newcommand{\der}[1]{\frac{\partial}{\partial #1}}
\newcommand{\norm}[1]{\parallel #1 \parallel}
\renewcommand{\matrix}[4]{\left( \begin{array}{cc}#1 & #2 \\ #3 & #4 \end{array}
          \right)}
\renewcommand{\sp}[2]{\langle #1,#2\rangle}

\newtheorem{theorem}{Theorem}[section]
\newtheorem{conjecture}[theorem]{Conjecture}
\newtheorem{lemma}[theorem]{Lemma}
\newtheorem{corollary}[theorem]{Corollary}
\newtheorem{proposition}[theorem]{Proposition}

\newpage
\section*{Introduction}
In the paper \cite{padgeom} the author has extended the theory of Selberg-type zeta functions to higher rank $p$-adic groups.
This extension remained incomplete insofar as within a higher rank group only elements of splitrank one were considered.\\
In the analogous setting of real Lie groups it emerged in recent years that the role of the Selberg zeta function in higher rank spaces is played by certain Lefschetz formulae attached to torus actions \cite{primgeo}.

We explain this in more detail.
Recall for rank one groups the relation between the Selberg zeta function and the trace formula.
It is clear that the analytical properties of the zeta function are derived by means of the trace formula.
Less well known is the fact that one can deduce the trace formula from the location of poles and zeros of the zeta function by evaluating a contour integral.
Likewise, in the case of higher rank groups and splitrank one elements, the Selberg zeta function corresponds to a Lefschetz formula attached to the action of a minimal split torus.
There is no proper analogue of the zeta function for higher rank elements, but there is a Lefschetz formula for that case, too.
A special version of the Lefschetz formula in the real setting was shown in \cite{primgeo} and a general version in \cite{Lefschetz}.
In the present paper we give the general Lefschetz formula in the $p$-adic setting.

\section{The trace formula}
Let $F$ be a nonarchimedean local field with valuation
ring $\CO$ and uniformizer $\varpi$. Denote by $G$ be a
semisimple linear algebraic group over $F$. Let $K\subset G$ be a
good maximal compact subgroup. Choose a parabolic
subgroup $P=LN$ of $G$ with Levi component $L$. 
Let $A$
denote the largest split torus in the center of $L$. 
Then
$A$ is called the {\it split component} of $P$. Let
$\Phi=\Phi(G,A)$ be the root system of the pair $(G,A)$,
i.e. $\Phi$ consists of all homomorphisms $\alpha :A\ra
\GL_1$ such that there is $X$ in the Lie algebra of $G$
with $\Ad(a)X= a^\alpha X$ for every $a\in A$. Given
$\alpha$, let $\n_\alpha$ be the Lie algebra generated by
all such $X$ and let $N_\alpha$ be the closed subgroup of
$N$ corresponding to $\n_\alpha$. Let $\Phi^+=\Phi(P,A)$ be
the subset of $\Phi$ consisting of all positive roots with
respect to $P$. Let $\Delta\subset\Phi^+$ be the subset of
simple roots. 
 Let $A^-\subset A$ be the set of all
$a\in A$ such that $|a^\alpha| <1$ for any $\alpha\in
\Delta$.

There is a reductive subgroup $M$ of $L$ with compact center such
that $MA$ has finite index in $L$. 
We can choose $M$ such that $K_M=M\cap K$ is a good
maximal subgroup of $M$. 
An element $g$ of $G$ is called \emph{elliptic} if it is
contained in a compact torus. 
Let $M_{ell}$ denote the set of elliptic elements
in $M$.

Let $X^*(A)=\Hom(A,\GL_1)$ be the group of all
homomorphisms as algebraic groups from $A$ to $\GL_1$. This
group is isomorphic to $\Z^r$ with $r=\dim A$. Likewise let
$X_*(A)=\Hom(\GL_1,A)$. There is a natural $\Z$-valued
pairing
\eqn
X^*(A)\times X_*(A) &\ra& \Hom(\GL_1,\GL_1)\cong\Z\\
(\alpha,\eta) &\mapsto& \alpha\circ\eta.
\endeqn
For every root $\alpha\in\Phi(A,G)\subset X^*(A)$ let
$\breve{\alpha}\in X_*(A)$ be its coroot. Then
$(\alpha,\breve{\alpha})=2$. The valuation $v$ of $F$
gives a group homomorphism $\GL_1(F)\ra\Z$. Let $A_c$ be
the unique maximal compact subgroup of $A$. Let $\Sigma
=A/A_c$; then $\Sigma$ is a $\Z$-lattice of rank $r=\dim
A$. By composing with the valuation $v$ the group $X^*(A)$
can be identified with
$$
\Sigma^*\=\Hom(\Sigma,\Z).
$$
Let
$$
\a_0^*\=\Hom(\Sigma,\R)\ \cong\ X^*(A)\otimes\R
$$
be the real vector space of all group homomorphisms from
$\Sigma$ to $\R$ and let
$\a^*=\a_0^*\otimes\C=\Hom(\Sigma,\C)\cong
X^*(A)\otimes\C$. For $a\in A$ and $\la\in\a^*$ let
$$
a^\la\=q^{-\la(a)},
$$
where $q$ is the number of elements in the residue class field of
$F$. In this way we get an identification
$$
{\a^*}/\mbox{\small$\frac{2\pi i}{\log q}$}\Sigma^* \
\cong\ \Hom(\Sigma,\C^\times).
$$
A quasicharacter $\nu : A\ra\C^\times$ is called {\it
unramified} if $\nu$ is trivial on $A_c$. The set
$\Hom(\Sigma,\C^\times)$ can be identified with the set of
unramified quasicharacters on $A$. Any unramified
quasicharacter $\nu$ can thus be given a unique real part
$\Re(\nu)\in \a_0^*$. This definition extends to not
necessarily unramified quasicharacters $\chi:A\ra\C^\times$
as follows. Choose a splitting $s:\Sigma\ra A$ of the exact
sequence
$$
1\ra A_c\ra A\ra\Sigma\ra 1.
$$
Then $\nu=\chi\circ s$ is an unramified character of $A$. Set
$$
\Re(\chi)\=\Re(\nu).
$$
This definition does not depend on the choice of the splitting
$s$. For quasicharacters $\chi$, $\chi'$ and $a\in A$ we will
frequently write $a^\chi$ instead of $\chi(a)$ and
$a^{\chi+\chi'}$ instead of $\chi(a)\chi'(a)$. Note that the
absolute value satisfies $|a^\chi|=a^{\Re(\chi)}$ and that a
quasicharacter $\chi$ actually is a character if and only if
$\Re(\chi)=0$.

Let $\Delta_P : P\ra\R_+$ be the modular function of the
group $P$. Then there is $\rho\in\a_0^*$ such that
$\Delta_P(a)=|a^{2\rho}|$. For $\nu\in\a^*$ and a root
$\alpha$ let
$$
\nu_\alpha\= (\nu,\breve{\alpha})\ \in\
X^*(\GL_1)\otimes\C\ \cong\ \C.
$$
Note that $\nu\in\a_0^*$ implies $\nu_\alpha\in\R$ for
every $\alpha$. For $\nu\in\a_0^*$ we say that $\nu$ is
positive, $\nu>0$, if $\nu_\alpha>0$ for every positive
root $\alpha$.

{\bf Example.} Let $G=\GL_n(F)$ and let $\varpi_j\in G$ be
the diagonal matrix
$\varpi_j=\diag(1,\dots,1,\varpi,1,\dots,1)$ with the
$\varpi$ on the $j$-th position. Let $\nu\in\a^*$ and let
$$
\nu_j\=\nu(\varpi_j A_c)\ \in\ \C.
$$
Let $\alpha$ be a root, say
$\alpha(\diag(a_1,\dots,a_n))=\frac{a_i}{a_j}$. Then
$$
\nu_\alpha\= \nu_i-\nu_j.
$$
Hence $\nu\in\a_0^*$ is positive if and only if
$\nu_1>\nu_2>\dots >\nu_n$.

We will fix Haar-measures of $G$ and its reductive subgroups as follows.
For $H\subset G$ being a torus there is a unique maximal compact subgroup $U_H$ which is open. 
Then we fix a Haar measure on $H$ such that $\vol (U_H)=1$.
If $H$ is connected reductive with compact center then we choose the unique positive Haar-measure which up to sign coincides with the Euler-Poincar\'e measure
\cite{Kottwitz}.
So in the latter case our measure is determined by the following property:
For any discrete torsionfree cocompact subgroup $\Ga_H\subset H$ we have
$$
\vol (\Ga_H \bs H) \= (-1)^{q(H)}\chi(\Ga_H ,\Q),
$$
where $q(H)$ is the $k$-rank of the derived group $H_{der}$ and $\chi(\Ga_H ,\Q)$ the
Euler-Poincar\'e characteristic of $H^\bullet(\Ga_H ,\Q)$. For the applications
recall that centralizers of tori in connected groups are connected
\cite{borel-lingroups}.

Assume we are given a discrete subgroup $\Ga$ of $G$ such that the quotient space
$\Ga \bs G$ is compact. Let $(\omega ,V_\omega)$ be a finite dimensional unitary
representation of $\Ga$ and let $L^2(\Ga \bs G,\omega)$ be the Hilbert space
consisting of all measurable functions $f: G \ra V_\omega$ such that $f(\ga x) =
\omega(\ga) f(x)$ and $|f|$ is square integrable over $\Ga \bs G$ (modulo null
functions). Let $R$ denote the unitary representation of $G$ on $L^2(\Ga \bs
G,\omega)$ defined by right shifts, i.e. $R(g) \ph (x) = \ph (xg)$ for $\ph \in
L^2(\Ga \bs G,\omega)$. It is known that as a $G$-representation this space splits
as a topological direct sum:
$$
L^2(\Ga \bs G,\omega) \= \bigoplus_{\pi \in \hat{G}} N_{\Ga ,\omega}(\pi) \pi
$$
with finite multiplicities $N_{\Ga ,\omega}(\pi)<\infty$.

Let $f$ be integrable over $G$, so $f$ is in $L^1(G)$.
The integral
$$
R(f) := \int_G f(x) R(x) \ dx
$$
defines an operator on the Hilbert space $L^2(\Ga \bs G,\omega)$.

For $g\in G$ and $f$ any function on $G$ we define the \emph{orbital integral}
$$
\CO_g(f) := \int_{G_g\bs G} f(x^{-1} gx)\ dx,
$$
whenever the integral exists.
Here $G_g$ is the centralizer of $g$ in $G$.
It is known that the group $G_g$ is unimodular, so we have an invariant measure on
$G_g\bs G$.

A function $f$ on $G$ or any of its closed subgroups is called
\emph{smooth} if it is locally constant. It is called \emph{uniformly smooth} if
there is an open subgroup
$U$ of
$G$ such that $f$ factors over
$U\bs G/U$. This is in particular the case if $f$ is smooth and compactly supported.

\begin{proposition}
(Trace formula)
Let $f$ be integrable and uniformly smooth, then we have 
$$
\sum_{\pi \in\hat{G}} N_{\Ga ,\omega}(\pi)\ \tr \pi (f) \= \sum_{[\ga]} \tr \omega (\ga)\ \vol(\Ga_\ga \bs G_\ga)\ \CO_\ga (f),
$$
where the sum on the right hand side runs over the set of $\Ga$-conjugacy classes $[\ga]$ in $\Ga$ and $\Ga_\ga$ denotes the centralizer of $\ga$ in $\Ga$.
Both sides converge absolutely and the left hand side actually is a finite sum.
\end{proposition}

\prf
At first fix a fundamental domain $\CF$ for $\Ga \bs G$ and let $\ph\in L^2(\Ga \bs G,\omega)$, then
\begin{eqnarray*}
R(f) &\=& \int_Gf(y) \ph(xy)\ dy\\
	&\=& \int_G f(x^{-1}y) \ph(y)\ dy\\
	&\=& \sum_{\ga\in\Ga} \int_\CF f(x^{-1}\ga y)\ph(\ga y)\ dy\\
	&\=& \int_{\Ga\bs G} \left( \sum_{\ga\in\Ga} f(x^{-1}\ga y)\omega (\ga)\right) \ph (y) \ dy
\end{eqnarray*}

We want to show that the sum $\sum_{\ga\in\Ga} f(x^{-1}\ga y)\omega (\ga)$ converges in $\End(V_\omega)$ absolutely and uniformly in $x$ and $y$.
Since $y$ can be replaced by $\ga y$, $\ga\in\Ga$ and since $\omega$ is unitary, we
only have to show the convergence  of $\sum_{\ga\in\Ga}|f(x^{-1}\ga y)|$ locally
uniformly in $y$. Let $\ga$ and $\tau$ be in $\Ga$  and assume that $x^{-1}\ga y$
and $x^{-1}\tau y$ lie in the same class in $G/U$. Then it follows $\tau yU\cap\ga
yU \ne \emptyset$ so with $V=yUy^{-1}$ we have $\ga^{-1} \tau V \cap V \ne
\emptyset$. It is clear that $V$ depends on $y$ only up to $U$ so to show locally
uniform convergence in $y$ it suffices to fix $V$. Since $V$ is compact also $V^2 =
\{ vv' |v,v'\in V\}$ is compact and so $\Ga \cap V^2$ is finite. This implies that
there are only finitely many $\ga\in\Ga$ with $\ga V\cap V\ne \emptyset$. Hence the
map $\Ga \ra G/U$, $\ga \mapsto x^{-1}\ga yU$ is finite to one with fibers having
$\le n$ elements for some natural number $n$. For $y$ fixed modulo $U$ we get
\begin{eqnarray*}
\sum_{\ga\in\Ga} |f(x^{-1} \ga y)| &\ \le\ & n\int_{G/U} |f(x)|\ dx\\
	&\=& \frac{n}{\vol(U)} \parallel f\parallel^1.
\end{eqnarray*}
We have shown the uniform convergence of the sum
$$
k_f(x,y) \= \sum_{\ga\in\Ga} f(x^{-1} \ga y) \omega(\ga).
$$

Observe that $R(f)$ factors over $L^2(\Ga\bs G,\omega)^U = L^2(\Ga\bs G/U,\omega)$, which is finite dimensional since $\Ga\bs G/U$ is a finite set.
So $R(f)$ acts on a finite dimensional space and $k_f(x,y)$ is the matrix of this operator.
We infer that $R(f)$ is of trace class, its trace equals
$$
\sum_{\pi\in\hat{G}} N_{\Ga ,\omega}(\pi)\ \tr\pi(f),
$$
and the sum is finite.
Further, since $k_f(x,y)$ is the matrix of $R(f)$ this trace also equals
\begin{eqnarray*}
\int_{\Ga \bs G} \tr k_f(x,x)\ dx &\=& \sum_{\ga\in\Ga} \int_\CF f(x^{-1}\ga x) \ dx\ \tr\omega(\ga)\\
	&\=& \sum_{[\ga]} \sum_{\sigma \in \Ga_\ga \bs \Ga} \int_\CF f((\sigma x)^{-1} \ga (\sigma x))\ dx\ \tr\omega(\ga)\\
	&\=& \sum_{[\ga]}  \int_{\Ga_\ga \bs G} f(x^{-1} \ga x)\ dx\ \tr\omega(\ga)\\
	&\=& \sum_{[\ga]}  \vol(\Ga_\ga\bs G_\ga)\ \int_{G_\ga \bs G} f(x^{-1} \ga x)\ dx\
\tr\omega(\ga).
\end{eqnarray*}
\qed

\section{The covolume of a centralizer}
Suppose $\ga\in\Ga$ is $G$-conjugate to some $a_\ga m_\ga\in A^+M_{ell}$.
We want to compute the covolume
$$
\vol(\Ga_\ga \bs G_\ga) .
$$
An element $x$ of $G$ is called 
 \emph{neat} if for  every
representation $\rho : G\ra GL_n(F)$ of $G$ the matrix $\rho(x)$ has no eigenvalue
which is a root of unity different from $1$. 
A subset $A$ of $G$ is called neat if each element of it
is. Every arithmetic
$\Ga$ has a finite index subgroup which is neat \cite{borel}. 

\begin{lemma}
Let $x\in G$ be neat and semisimple. Let $G_x$ denote its centralizer in $G$. Then
for every
$k\in\N$ we have $G_x=G_{x^k}$.
\end{lemma}

\prf
Since $G$ is linear algebraic it is a subgroup of some $H=\GL_n(\bar F)$, where $\bar
F$ is an algebraic closure of $F$. If we can show the claim for $H$ then it follows
for $G$ as well since $G_x=H_x\cap G$. In $H$ we can assume $x$ to be a diagonal
matrix. Since
$x$ is neat this implies the claim.
\qed

We suppose that $\Ga$ is
neat. This implies that for any
$\ga\in\Ga$ the Zariski closure of the group generated by $\ga$ is a torus. It then
follows  that $G_\ga$ is a connected reductive group \cite{borel-lingroups}.

An element $\ga\in\Ga$ is called \emph{primitive} if $\ga =\sigma^n$ with
$\sigma\in\Ga$ and $n\in\N$ implies $n=1$. It is a property of discrete cocompact
torsion free subgroups $\Ga$ of $G$ that every $\ga\in\Ga$, $\ga\ne 1$ is a positive
power of a unique primitive element. In other words, given a nontrivial $\ga\in\Ga$
there exists a unique primitive $\ga_0$ and a unique $\mu(\ga)\in\N$ such that
$$
\ga =\ga_0^{\mu(\ga)}.
$$

Let $\Sigma$ be a group of finite cohomological dimension $cd(\Sigma)$ over $\Q$.
We write
$$
\chi(\Sigma) \= \chi(\Sigma ,\Q) \ :=\ \sum_{p=0}^{cd(\Sigma)} (-1)^p \dim H^p(\Sigma
,\Q),
$$
for the Euler-Poincar\'e characteristic.
We also define the higher Euler characteristic as
$$
\chi_{_r}(\Sigma) \= \chi_{_r}(\Sigma ,\Q) \ :=\ \sum_{p=0}^{cd(\Sigma)} 
(-1)^{p+r}\binom pr
\dim H^p(\Sigma ,\Q),
$$
for $r=1,2,3,\dots$
It is known that $\Ga$ has finite cohomological dimension over $\Q$.

 We denote by $\CE_P(\Ga)$
the set of all conjugacy classes $[\ga]$ in $\ga$ such that $\ga$ is in $G$
conjugate to an element $a_\ga m_\ga\in AM$, where $m_\ga$ is elliptic and
$a_\ga\in A^-$.

Let $\ga\in\CE_P(\Ga)$. To simplify the notation let's assume that $\ga=a_\ga
m_\ga\in A^- M_{ell}$. Let $C_\ga$ be the connected component of the center of
$G_\ga$ then
$C_\ga = AB_\ga$, where
$B_\ga$ is the connected center of $M_{m_\ga}$ the latter group will also be
written as $M_\ga$. Let $M_\ga^{der}$ be the derived group of $M_\ga$. Then
$M_\ga=M_\ga^{der} B_\ga$.

\begin{lemma}
$B_\ga$ is compact.
\end{lemma}

\prf
Since $m_\ga$ is elliptic there is a compact Cartan subgroup $T$ of $M$ containing
$m_\ga$. Since $M$ modulo its center is a connected semisimple linear algebraic
group it follows that $T$ is a torus and therefore abelian. Therefore $T\subset
M_{m_\ga}$. Let
$b\in B_\ga$. Then $b$ commutes with every $t\in T$, therefore $b$ 
lies in the centralizer of $T$ in $M$ which equals $T$. So we have shown
$B_\ga\subset T$.
\qed

Let $\Ga_{\ga,A}=A\cap \Ga_\ga B_\ga$ and $\Ga_{\ga,M}=M_\ga^{der}\cap \Ga_\ga
AB_\ga$. Similar to the proof of Lemma 3.3 of \cite{Wolf} one shows that
$\Ga_{\ga,A}$ and
$\Ga_{\ga,M}$ are discrete cocompact subgroups of $A$ and $M_\ga^{der}$ resp.
Let
$$
\la_{\ga}\df \vol(\Ga_{\ga,A}\bs A).
$$

\begin{proposition}\label{2.3}
Assume $\Ga$ neat and let $\ga\in\Ga$ be $G$-conjugate to an element of $A^+M_{ell}$.
Then we get
$$
\vol (\Ga_\ga\bs G_\ga) = \la_\ga\ (-1)^{q(G)+r}\ \chi_{_r}(\Ga_\ga),
$$
where $r=\dim A$.
\end{proposition}

\prf
We normalize the volume of $B_\ga$ to be $1$. Then
\begin{eqnarray*}
\vol(\Ga_\ga\bs G_\ga) &=& \vol(\Ga_\ga\bs AM_\ga)\\
&=& \vol(\Ga_\ga B_\ga \bs AM_\ga)
\end{eqnarray*}
The space $\Ga_\ga B_\ga \bs AM_\ga$ is the total space of a fibration with fibre
$\Ga_{\ga,A}\bs A$ and base space $\Ga_\ga AB_\ga\bs M_\ga A\cong \Ga_{\ga,M}\bs
M_\ga^{der}$. Hence
$$
\vol(\Ga_\ga B_\ga \bs AM_\ga) \= \vol (\Ga_{\ga,A}\bs A) \,\vol(\Ga_{\ga,M}\bs
M_\ga^{der}).
$$
Since $\la_\ga=\vol (\Ga_{\ga,A}\bs A)$ it remains to show
$$
\vol(\Ga_{\ga,M}\bs M_\ga^{der})\= (-1)^r\chi_r(\Ga_\ga).
$$
We know that
$$
\vol(\Ga_{\ga,M}\bs M_\ga^{der})\= (-1)^{q(M_\ga)}\chi(\Ga_{\ga,M})\=
(-1)^{q(G)+r}\chi(\Ga_{\ga,M}).
$$
So it remains to show that $\chi(\Ga_{\ga,M})=\chi_r(\Ga_\ga)$.
The group $\Ga_{\ga,M}$ is isomorphic to $\Ga_\ga/\Sigma$, where $\Sigma=\Ga\cap
AB_\ga$ is isomorphic to $\Z^r$. So the proposition follows from the next Lemma.

\begin{lemma}\label{chichi1}
Let $\Ga,\Lambda$ be of finite cohomological dimension over $\Q$.
Let $C_r$ be a group isomorphic to $\Z^r$ and assume there is an exact sequence
$$
1\ra C_r\ra \Ga\ra \Lambda\ra 1.
$$
Assume that $C_r$ is central in $\Ga$. 
Then 
$$
\chi (\Lambda,\Q) \= \chi_{r}(\Ga,\Q).
$$
\end{lemma}

\prf
We first consider the case $r=1$. In this case we want to prove for every $r$,
$$
\chi_{r-1} (\Lambda,\Q) \= \chi_{r}(\Ga,\Q).
$$
For this consider the Hochschild-Serre spectral sequence:
$$
E_2^{p,q} \= H^p(\Lambda,H^q(C_1,\Q))
$$
which abuts to
$$
H^{p+q}(\Ga,\Q).
$$
Since $C_1\cong \Z$ it follows
$$
H^q(C_1,\Q) \= \left\{ \begin{array}{cl} \Q& {\rm if}\ q=0,1\\ 0&{\rm
else}.\end{array}\right.
$$
Since $C_1$ is infinite cyclic and central it is an exercise to see that the
spectral sequence degenerates at $E_2$. Therefore,
\begin{eqnarray*}
\chi_r(\Ga) &=& \sum_{j\ge 0} (-1)^{j+r}\binom jr \dim H^j(\Ga)\\
&=& \sum_{j\ge r} (-1)^{j+r}\binom jr (\dim H^j(\Lambda)+\dim H^{j-1}(\Lambda)\\
&=& \sum_{j\ge r} (-1)^{j+r}\binom jr \dim H^j(\Lambda)\\
&& \quad -\sum_{j\ge{r-1}} (-1)^{j+r} \binom{j+1}r \dim H^j(\Lambda).
\end{eqnarray*}
Now replace $\binom{j+1}r$ by $\binom jr +\binom j{r-1}$ to get the claim. For the
general case write $C_r=C_1\oplus C^1$, where $C_1$ is cyclic and $C^1\cong
\Z^{r-1}$. Apply the above to $C_1$ and iterate this to get the lemma and hence the
proposition.
\qed

\section{The Lefschetz formula}
For a representation $\pi$ of $G$ let $\pi^\infty$ denote the subrepresentation of
\emph{smooth vectors}, ie $\pi^\infty$ is the representation on the space
$\bigcup_{H\subset G} \pi^H$, where $H$ ranges over the set of all open subgroups
of $G$. Further let $\pi_N$ denote the \emph{Jacquet module} of $\pi$. By
definition $\pi_N$ is the largest quotient $MAN$-module of $\pi^\infty$ on which
$N$ acts trivially. One can achieve this by factoring out the vector subspace
consisting of all vectors of the form $v-\pi(n)v$ for $v\in\pi^\infty$, $n\in N$. It
is known that if $\pi$ is an irreducible admissible representation, then $\pi_N$ is a
admissible $MA$-module of finite length. For a smooth $M$-module $V$ let
$H_c^\bullet(M,V)$ denote the continuous cohomology with coefficients in $V$ as in
\cite{Borel-Wallach}.

\begin{theorem}(Lefschetz Formula)\\
Let $\Ga$ be a neat discrete cocompact subgroup of $G$.
Let $\ph$ be a uniformly smooth function on $A$ with support in $A^-$. Suppose that
the function $a\mapsto \ph(a)|a^{-2\rho}|$ is integrable on $A$. Let
$\sigma$ be a finite dimensional representation of $M$. Let $q$ be the $F$-splitrank
of $G$ and $r=\dim A$.Then
$$
\sum_{\pi\in\hat G} N_{\Ga,\omega}(\pi)\sum_{q=0}^{\dim M} (-1)^a \int_{A^-}
\ph(a)\,\tr(a| H_c^q(M,\pi_N\otimes\sigma))\, da
$$
equals
$$
(-1)^{q+r}\sum_{[\ga]\in\CE_P(\Ga)} \la_\ga\,
\chi_r(\Ga_\ga)\,\tr\omega(\ga)\,\tr\sigma(m_\ga)\,\ph(a_\ga)\,|a_\ga^{2\rho}|.
$$
Both outer sums converge absolutely and the sum over $\pi\in\hat G$ actually is a
finite sum, ie, the summand is zero for all but finitely many $\pi$. For a given
compact open subgroup $U$ of $A$ both sides represent a continuous linear functional
on the space of all functions $\ph$ as above which factor over $A/U$, where this
space is equipped with the norm $\norm \ph=\int_A|\ph(a)| |a^{-2\rho}|\, da$.
\end{theorem}

Let $A^*$ denote the set of all continuous group homomorphisms $\la\colon
A\ra\C^\times$. For $\la\in A^*$ and an $A$-module $V$ let $V_\la$ denote the
generalized $\la$-eigenspace, ie,
$$
V_\la\df \bigcup_{k=1}^\infty \{ v\in V\mid (a-\la(a))^k v=0\ \forall a\in A\}.
$$
Then
$$
\int_{A^-} \ph(a)\,\tr(a|H_c^q(M,\pi_N\otimes\sigma))\, da\=\sum_{\la\in A^*}\dim
H_c^q(M,\pi_N\otimes\sigma)_\la\,\int_{A^-}\ph(a\,\la(a\, da.
$$
For $\la\in A^*$ define
$$
m_\la^{\sigma,\omega} \df \sum_{\pi\in\hat G}N_{\Ga,\omega}(\pi)\sum_{q=0}^{\dim
M}(-1)^q\,\dim H_c^q(M,\pi_N\otimes\sigma)_\la.
$$
The sum is always finite.

On the other hand, for $[\ga]\in\CE_P(\Ga)$ let
$$
c_\ga\df \la_\ga\,\chi_r(\Ga_\ga)\,|a_\ga^{2\rho}|.
$$
Then the Theorem is equivalent to the following Corollary.

\begin{corollary}
(Lefschetz Formula)\\
As an identity of distributions on $A^-$ we have
$$
\sum_{\la\in A^*} m_\la^{\sigma,\omega}\, \la\= \sum_{[\ga]\in\CE_P(\Ga)}
c_\ga\,\tr\omega(\ga)\,\tr\sigma(m_\ga)\,\delta_{a_\ga}.
$$
\end{corollary}

{\bf Proof of the Theorem:}
Let $f_{EP}$ be an Euler-Poincar\'e function on $M$ which
is $K_M$-central \cite{Kottwitz}. For $m\in M$ regular we have
$$
\CO_m^M(f_{EP})\=\begin{cases} 1 & m {\ \rm elliptic},\\ 0 & {\rm
otherwise}.\end{cases}
$$
For $g\in G$ and a finite dimensional $F$ vector space  $V$ on which $g$ acts
linearly let
$E(g|V)$ be the set of all absolute values $|\mu|$, where $\mu$ ranges over
the eigenvalues of $g$ in the algebraic closure $\bar F$ of $F$. 
Let $\la_{min}(g|V)$ denote the minimum and $\la_{max}(v|V)$ the maximum of
$E(g|V)$. For $am\in AM$ define
$$
\la(am)\df \frac{\la_{min}(a|\bar\n)}{\la_{max}(m|\g)^2}
$$
Note that $\la_{max}(m|\g)$ is always $\ge 1$ and that
$\la_{max}(m|\g)\la_{min}(m|\g)=1$. We will consider the set
$$
(AM)^\sim \ :=\ \{ am\in AM | \la(am)>1 \}.
$$
Let $M_{ell}$ denote the set of elliptic elements in $M$.

\begin{lemma} \label{MA}
The set $(AM)^\sim$ has the following properties:

\begin{enumerate}
\item
$A^-M_{ell}\subset (AM)^{\sim}$
\item
$am\in (AM)^\sim \Rightarrow a\in A^-$
\item
$am, a'm' \in (AM)^\sim, g\in G\ {\rm with}\ a'm'=gamg^{-1}
\Rightarrow a=a', g\in AM$.
\end{enumerate}
\end{lemma}

\prf The first two are immediate. For the third let $am, a'm' \in
(AM)^\sim$ and $g\in G$ with $a'm'=gamg^{-1}$. Observe that by the
definition of $(AM)^\sim$ we have
\begin{eqnarray*}
\la_{min}(am | \bar{\n}) &\ge& \la_{min}(a| \bar{\n})\la_{min}(m |
\g)\\
 &>& \la_{max}(m|\g)^2\la_{min}(m|\g)\\
 &=& \la_{max}(m|\g)\\
 &\ge & \la_{max}(m | \a +\m +\n)\\ &\ge&\la_{max}(am | \a +\m +\n)
\end{eqnarray*}
that is, any eigenvalue of $am$ on $\bar{\n}$ is strictly bigger
than any eigenvalue on $\a +\m +\n$. Since $\g = \a +\m +\n
+\bar{\n}$ and the same holds for $a'm'$, which has the same
eigenvalues as $am$, we infer that $\Ad(g)\bar{\n} =\bar{\n}$. So
$g$ lies in the normalizer of $\bar{\n}$, which is
$\bar{P}=MA\bar{N} =\bar{N}AM$. Now suppose $g=nm_1a_1$ and
$\hat{m} =m_1mm_1^{-1}$ then
$$
gamg^{-1} \= na\hat{m}n^{-1} \= a\hat{m}\
(a\hat{m})^{-1}n(a\hat{m})\ n^{-1}.
$$
Since this lies in $AM$ we have $(a\hat{m})^{-1}n(a\hat{m}) =n$
which since $am\in (AM)^\sim$ implies $n=1$. The lemma is proven.
\qed

 Let $G$ act on itself
by conjugation, write $g.x = gxg^{-1}$, write $G.x$ for the orbit,
so $G.x = \{ gxg^{-1} | g\in G \}$ as well as $G.S = \{ gsg^{-1} |
s\in S , g\in G \}$ for any subset $S$ of $G$.

Fix a smooth function $\eta$ on $N$ which has compact support, is
positive, invariant under $K_M$ and satisfies $\int_N\eta(n) dn
=1$. Extend the function $\ph$ from $A^-$ to a conjugation
invariant smooth function $\tilde{\ph}$ on $AM$ such that
$\tilde{\ph}(am)=\ph(a)$ whenever $m$ is elliptic and such that
there is a compact subset $C\subset A^-$ such that $\tilde{\ph}$
is supported in $CM\cap(AM)^\sim$.
 It follows that the function
$$
am\ \mapsto\ f_{EP}(m)\,\tr\sigma(m)\,\tilde{\ph}(am)\,|a^{2\rho}|
$$
is smooth and integrable on $AM$. Given these data let
$f = f_{\eta ,\tau ,\ph} : H\ra \C$ be defined by
$$
f (kn ma (kn)^{-1}) := \eta (n) f_{EP}(m)
\,\tr\sigma(m)\,\tilde{\ph}(am)\,|a^{2\rho}|,
$$
for $k\in K, n\in N, m\in M, a\in\overline{A^-}$. Further $f(x)=0$
if $x$ is not in $G.(AM)^\sim$.

\begin{lemma} \label{welldef}
The function $f$ is well defined.
\end{lemma}

\prf By the decomposition $G=KP=KNMA$ every element $x\in
G.(AM)^\sim$ can be written in the form $kn ma (kn)^{-1}$. Now
suppose two such representations coincide, that is
$$
kn ma (kn)^{-1}\= k'n' m'a' (k'n')^{-1}
$$
then by Lemma \ref{MA} we get $(n')^{-1} (k')^{-1}kn\in MA$, or
$(k')^{-1}k\in n'MAn^{-1}\subset MAN$, hence $(k')^{-1}k\in K\cap
MAN=K\cap M=K_M$. Write $(k')^{-1}k=k_M$ and $n''=k_Mnk_M^{-1}$,
then it follows
$$
n'' k_Mmk_M^{-1} a (n'')^{-1}\= n'm'a'(n')^{-1}.
$$
Again by Lemma \ref{MA} we conclude $(n')^{-1}n''\in MA$, hence
$n'=n''$ and so
$$
k_Mmk_M^{-1} a \= m'a',
$$
which implies the well-definedness of $f$.
\qed

We will plug $f$ into the trace formula. For the geometric side
let $\ga \in \Ga$. We have to calculate the orbital integral:
$$
\CO_\ga (f) = \int_{G_\ga \bs G} f(x^{-1}\ga x) dx.
$$
by the definition of $f$ it follows that $\CO_\ga(f)=0$ if
$\ga\notin G.(AM)^\sim$. It remains to compute $\CO_{am}(f)$ for
$am\in(AM)^\sim$. Again by the definition of $f$ it follows
\begin{eqnarray*}
\CO_{am}(f) &=&
\CO_m^M(f_{EP})\,\tr\sigma(m)\,\tilde{\ph}(am)\,|a^{2\rho}|\\
&=& \begin{cases} \tr\sigma(m)\,\ph(a)\,|a^{2\rho}| & {\rm if}\ m\ {\rm is\
elliptic},\\ 0 & {\rm otherwise.}\end{cases}
\end{eqnarray*}
Here $\CO_m^M$ denotes the orbital integral in the group $M$.
 Recall that Proposition \ref{2.3} says
$$
\vol (\Ga_\ga\bs G_\ga) = (-1)^{q(G)+r}\,\la_\ga\,\chi_{_r}(\Ga_\ga),
$$
so that for $\ga\in\CE_P(\Ga)$,
$$
\vol(\Ga_\ga\bs G_\ga)\,\CO_\ga(f)\=
(-1)^{q(G+r}\,\la_\ga\,\chi_r(\Ga_\ga)\,\tr\sigma(m_\ga)\,\ph(a_\ga)\,|a_\ga^{2\rho}|.
$$

To compute the spectral side let $\pi\in\hat{G}$. We want to compute
$\tr\pi(f)$. Let $\Theta_\pi^G$ be the locally integrable
conjugation invariant function  on $G$ such that
$$
\tr\pi(f)\= \int_G f(x)\, \Theta_\pi^G(x) dx.
$$
This function $\Theta_\pi$ is called the \emph{character} of $\pi$.
It is known that the Jacquet module $\pi_{{N}}$ is a finitely generated admissible
module for the  group $MA$ and therefore it has a character
$\Theta_{\pi_{{N}}}^{MA}$. In \cite{Casselman} it is shown that
$$
\Theta_\pi(am) \= \Theta_{\pi_{{N}}}^{MA}(ma)
$$
for $ma\in A^-M_{ell}$.

Let $h$ be a function in $L^1(G)$ which is supported in the set $G.MA$. Comparing
invariant differential forms as in the proof of the Weyl integration formula one
gets that the integral $\int_G h(x)\, dx$ equals
$$
\frac 1{|W(G,A)|}
\int_A\int_M\int_{G/AM} h(y amy^{-1})\,|\det(1-am | \n +\bar \n)| dy\, da\, dm,
$$
where $W(G,A)$ is the Weyl group of $A$ in $G$. 

For $a\in A^-$ and $m\in M_{ell}$ every eigenvalue of $am$ on $\n$ is of absolute
value $<1$ and $>1$ on $\bar \n$. By the ultrametric ptoperty this implies
\begin{eqnarray*}
|\det(1-am|\n +\bar \n)| &=& |\det(1-am|\bar \n)|\\
&=& |\det(am|\bar \n)|\\
&=& |\det(a|\bar \n)|\\
&=& |a^{-2\rho}|,
\end{eqnarray*}

We apply this to
$h(x)=f(x)\Theta_\pi^G(x)$ and use conjugation invariance of $\Theta_\pi^G$ to get
that $\tr\pi(f)$ equals
$$
\frac 1{|W(G,A)|}\int_{AM}
f_{EP}(m)\,\tr\sigma(m)\,\tilde\ph(am)\,\Theta_{\pi_N}^M(am)\,da\, dm,
$$
which is the same as
$$
\int_{A^-M}
f_{EP}(m)\,\tr\sigma(m)\,\tilde\ph(am)\,\Theta_{\pi_N}^M(am)\,da\, dm.
$$
We recall the Weyl integration formula for $M$. Let $(H_j)_j$ be a maximal family of
pairwise non-conjugate Cartan subgroups of $M$. Let $W_j$ be the Weyl group of $H_j$
in $M$. For $h\in H_j$ let $D_j(h)=\det(1-h|\m/\h_j)$, where $\m$ and $\h_j$ are the
Lie algebras of $M$ and $H_j$ resp. Then, for every $h\in L^1(M)$,
\begin{eqnarray*}
\int_M h(m)\, dm&=& \sum_j \frac 1{|W_j|}\int_{H_j^{reg}}\int_{M/H_j}h(mxm^{-1})\,
D_j(x)\, dm\, dx\\
&=& \sum_j \frac 1{|W_j|}\int_{H_j} \CO_x^M(h)\, dx,
\end{eqnarray*}
where $H_j^{reg}$ is the set of $x\in H_j$ which are regular in $M$. We fix
$a\in A^-$ and apply this to
$h(m)=f_{EP}(m)\tilde\ph(am)\tr\sigma(m)\Theta_{\pi_N}^M(am)$.
Since
$\tilde\ph$ is conjugation invariant we get for $x\in H_j^{reg}$,
$$
\CO_x^M(h)\=\CO_x^M(f_{EP})\,\tilde\ph(am)\,\tr\sigma(m)\,\Theta_{\pi_N}^{AM}(ax).
$$
This is non-zero only if $x$ is elliptic. If $x$ is elliptic, then $\tilde\ph(ax)$
equals $\ph(a)$. So we can replace $\tilde\ph(ax)$ by $\ph(a)$ throughout. Thus
$\tr\pi(f)$ equals
$$
\int_{A^-M}
f_{EP}(m)\,\ph(a)\,\tr\sigma(m)\,\Theta_{\pi_N}^M(am)\,da\, dm.
$$
 The trace
$\tr\pi(f)$ therefore equals
$$
\int_{A^-M} f_{EP}(m)\, \ph(a)\,\Theta_{\pi_N\otimes\sigma}(am)\, da\,
dm.
$$
We write $H_c^\bullet(M,V)$ for the continuous cohomology of $M$ with coefficients
in the $M$-module $V$. By Theorem 2 in \cite{Kottwitz},
$$
\tr(\pi_N\otimes\sigma)(f_{EP})\=\sum_{q=0}^{\dim M} (-1)^q\, \dim
H_c^q(M,\pi_N\otimes\sigma).
$$
The cohomology groups $H_c^q(M,\pi_N\otimes\sigma)$ are finite dimensional
$A$-modules and
$$
\tr\pi(f)\=\sum_{q=0}^{\dim M} (-1)^q\,
\int_{A^-}\tr(a|H_c^q(M,\pi_N\otimes\sigma))\,\ph(a)\, da.
$$
The Lefschetz Theorem follows.
\qed

{\small
University of Exeter, Mathematics, Exeter EX4 4QE, Devon, UK\\
a.h.j.deitmar@ex.ac.uk}

\end{document}